\documentclass[11pt]{article}

\usepackage{natbib}
\usepackage{amsthm, amsmath, amssymb, epsfig}   
\usepackage{algorithm,algorithmic}
\usepackage[justification=centerfirst, singlelinecheck=true]{caption}

\newcommand{\oa}{\overline{a}}
\newcommand{\oaa}{\overline{\alpha}}
\newcommand{\ox}{\overline{x}}
\newcommand{\gb}{\operatorname{gb}}
\newcommand{\lpp}{\operatorname{lpp}}

\newcommand{\lex}{\operatorname{lex}}
\newcommand{\grevlex}{\operatorname{grevlex}}

\newcommand{\V}{\mathbb V}

\newcommand{\Compl}{\mathbb C}

\newcommand{\Real}{\mathbb R}
\newcommand{\anul}[1]{}

\newtheorem{theorem}{Theorem}

\theoremstyle{definition}
\newtheorem{definition}[theorem]{Definition}
\newtheorem{example}[theorem]{Example}

\title{Minimal Canonical Comprehensive Gr\"obner Systems\footnote{This research was partly supported by the
    Ministerio de Ciencia y Tecnolog\'{\i}a under project
    MTM 2006-01267, and by the Generalitat de
    Catalunya under project 2005 SGR 00692.}
}

\author{Montserrat Manubens \& Antonio Montes\\ \\
Departament de Matem\`atica Aplicada 2,\\ Universitat
Polit\`ecnica de Catalunya, Spain.\\
e-mail: \{montserrat.manubens, antonio.montes\}@upc.edu\\
http://www-ma2.upc.edu/$\sim$montes}

\begin{document}
\maketitle

\begin{abstract}
This is the continuation of Montes' paper ``On the canonical
discussion of polynomial systems with parameters". In this paper
we define the Minimal Canonical Comprehensive Gr\"obner System
(MCCGS) of a parametric ideal and fix under which hypothesis it
exists and is computable. An algorithm to obtain a canonical
description of the segments of the MCCGS is given, completing so
the whole MCCGS algorithm (implemented in Maple). We show its high
utility for applications, like automatic theorem proving and
discovering, and compare it with other existing methods. A way to
detect a counterexample is outlined, although the high number of
tests done give evidence of the existence of the MCCGS.



\end{abstract}

\section{Introduction}\label{intro}
In this paper we continue the task introduced in~\cite{Mo07}. Let
us briefly remember the basic features.

Given a parametric polynomial ideal $I\subset K[\oa][\ox]$ in the
variables $\ox=(x_1,\dots,x_n)$ and the parameters
$\oa=(a_1,\dots,a_m)$, and monomial order $\succ_{\ox}$, our
interest is to find the different types of solutions for the
different values of the parameters. Let $K$ be a computable field
and $\overline{K}$ an algebraically closed extension. A
specialization is the homomorphism
$\sigma_{\oaa}:K[\oa][\ox]\rightarrow {{\overline{K}}}[\ox]$, that
corresponds to the substitution of the parameters by concrete
values $\oaa \in {\overline{K}}^m$.  A comprehensive Gr\"obner
system (CGS) is a set of pairs:
\[
\begin{array}{lcl}
\operatorname{CGS}(I,\succ_{\ox})& = &
 \{ (S_i,B_i)\ : S_i\subseteq \overline{K}^m \text{constructible sets}, \ B_i
\subset A[\ox],\\
&& \ \sigma_{\oaa}(B_i) = \gb( \sigma_{\oaa}(I), \succ_{\ox}) \
\forall \oaa \in S_i\ ,\text{and } \bigcup_i S_i=\overline{K}^m
\},
\end{array}
\]
where the $S_i$ are called ``segments" and the $B_i$ ``bases".
Frequently the word ``segment" is also used for the pair
$(B_i,S_i)$ whenever the sense is clear from the context.

There are different known algorithms that provide Comprehensive
Gr\"obner Bases and Systems for a given ideal:
\begin{itemize}
  \item [] CGB~\cite{We92},
  \item [] ACGB~\cite{SaSu03,Sa05,Na05},
  \item [] SACGB~\cite{SuSa06},
  \item [] HSGB~\cite{GoTrZa05},
  \item [] BUILDTREE~\cite{Mo02, MaMo06,Mo07}.
\end{itemize}
There are available implementations~\citep{DoSeSt06} of
Weispfenning's CGB algorithm in Reduce, of Suzuki-Sato's SACGB in
Risa/Asir and in Maple~\cite{SuSa06} and of Montes's BUILDTREE in
Maple\footnote{The library DPGB 7.0 written in Maple 8 is
available at the web http://www-ma2.upc.edu/$\sim$montes, and is
actualized with the MCCGS algorithm.}. All these algorithms allow
to build both Comprehensive Gr\"obner Bases and Systems, but they
are differently oriented. A comparison of the most interesting
among them is given in section~\ref{compare}.

In fact, comprehensive Gr\"obner systems are in general more
effective to handle for their use in the applications than
comprehensive Gr\"obner bases. But it is also convenient to
require some more additional features to these Gr\"obner systems
when looking for applications.

The first requirement is to have {\em disjoint} and {\em reduced}
CGS. By disjoint we mean that the $S_i$ form a partition of
${\overline{K}}^m$, and by reduced that the bases $B_i$ specialize
to the reduced Gr\"obner basis of $\sigma_{\oaa}(I)$ preserving
the leading power products ($\lpp$), for every value $\oaa$ of the
parameters inside $S_i$. The algorithm BUILDTREE (introduced
in~\cite{Mo02} as DISPGB and improved in~\cite{MaMo06}) already
builds a disjoint, reduced CGS.

In~\cite{Mo07} the interest is focused in the improvement of
BUILDTREE to obtain a simpler and canonical CGS. The method
consists of grouping together all the segments with the same
$\lpp$ that allow a same basis specializing well on all the
grouped segments. A natural conjecture establishes the existence
of an equivalence relation between the segments having the same
$\lpp$, and an algorithm is given to compute the basis
corresponding to the grouped segments.

In order to obtain a truly canonical CGS we need to describe the
segments in a canonical way. This is the objective of the present
paper. In~\cite{Mo07} a canonical description of a segment
determined by a diff-specification was already given, but it
remained to obtain a canonical representation of the addition of
such segments. The objective is thus to obtain the MCCGS (minimal
canonical CGS).
\begin{definition}\label{DefMCCGS} We call Minimal Canonical CGS (MCCGS)
a CGS with the following properties:
\begin{enumerate}\renewcommand*{\labelenumi}{\roman{enumi})}
  \item disjoint CGS, i.e. $S_i \cap S_j =\emptyset$ for
  $i\not=j$;
  \item reduced CGS, i.e. the polynomials in $B_i$ have
  content $1$ w.r.t. $\ox$, $B_i$ specializes to the reduced Gr\"obner basis of
  $\sigma_{\oaa}(I)$ for every $\oaa \in S_i$, their leading coefficients are non-null
  on $S_i$ and their $\lpp$ remain stable;
  \item  the sets $S_1,\dots,S_s$ are intrinsic for the given $I$ and
  $\succ_{\ox}$ and are described in a canonical form.
  \item the number of segments of the CGS with the above properties is minimal.
\end{enumerate}
\end{definition}
The currently existing algorithms that can build comprehensive
Gr\"obner systems, say BUILDTREE, CGB, CCGB, ACGB and SACGB, do
not hold all these properties. BUILDTREE builds a comprehensive
Gr\"obner system satisfying properties i) and ii). But CGB, ACGB
and SACGB do not hold property i), at least. Finally, although the
Gr\"obner system obtained within CCGB is canonically determined,
does not hold properties i) nor ii) as for the obtention of a
comprehensive Gr\"obner basis the algorithm needs the Gr\"obner
systems to be faithful.

It must be emphasized that the existence of the MCCGS depends on
the Conjecture formulated in~\cite{Mo07} about the existence of an
equivalence relation between segments allowing a common basis.

If the Conjecture is true, then the computation using MCCGS
algorithm proposed in~\cite{Mo07} and in this paper, already
depends on the semi-algorithm GENIMAGE given there for computing
pre-images, that uses arbitrary bounds.

With these restrictions, MCCGS algorithm builds a comprehensive
Gr\"obner system satisfying all the properties in
Definition~\ref{DefMCCGS}. These properties will make the
algorithm more suitable for the applications. In particular, they
are very appropriate for automatic theorem proving and discovery
(see \cite{MoRe07}) as well as to compute geometric loci as shown
in example~\ref{conicex}.

Furthermore, MCCGS also allows to restrict the parameter space to
a constructible set and impose a-priori null and non-null
conditions. This is also an interesting tool for applications when
we want some degenerate cases to be avoided (see
Section~\ref{applicationsection}) or restrictions on the
parameters to be given. For example, when the parameters involve
angles, and the equations are given using the sine and cosine of
the angles as parameters, it is important to restrict the
solutions to $\cos^2 \varphi + \sin^2 \varphi -1=0$.

The whole algorithm MCCGS is achieved by three steps:
\begin{enumerate}\renewcommand*{\labelenumi}{\roman{enumi})}
  \item BUILDTREE (described in~\cite{MaMo06}),
  \item grouping segments with common basis (described
  in~\cite{Mo07}),
  \item representing the subsets in canonical form. This part will be
  described in sections~\ref{addingsection} and~\ref{algorithmsection}.
\end{enumerate}
Although the algorithm requires two term orders (one for the
variables $\succ_{\ox}$ and another for the parameters
$\succ_{\oa}$), the result will not depend on $\succ_{\oa}$, as
the segments $(B_i,S_i)$ are intrinsic for the given ideal $I$ and
the term order $\succ_{\ox}$. Even though, $\succ_{\oa}$ will be
used to determine the reduced Gr\"obner bases of the ideals
involved in the description of $S_i$.

The paper is structured as follows: section \ref{SecPreliminaries}
is devoted to recalling some properties and results from
\cite{Mo07} which are used in the subsequent sections. The
generalization of the canonical specification and its properties
are given in Section \ref{addingsection}. In Section
\ref{algorithmsection} we give the algorithm which collects the
corresponding segments into a generalized canonical specification
and builds up the Minimal Canonical Comprehensive Gr\"obner System
(MCCGS). In section \ref{applicationsection} a practical
application to automatic theorem proving is given. Finally, in
section \ref{compare} we compare the main available CGS
algorithms.

\section{Preliminaries}\label{SecPreliminaries}
We describe now briefly steps i) and ii) of the MCCGS before
tackling the last step iii) that will be studied in this paper.
The algorithm starts with a parametric ideal $I$ and a term-order
$\succ_{\ox}$ on the variables. An auxiliary term-order
$\succ_{\oa}$ over the parameters is needed to describe the
subsets in ${\overline{K}}^m$ using Gr\"obner bases. It does not
affect the segments themselves but only their description.

Step i) is performed by BUILDTREE algorithm, and was described for
the first time in~\cite{Mo02} and improved in~\cite{MaMo06}. The
output is a disjoint reduced CGS, where the subsets $S_i$ are
determined by {\em red-specifications}. A red-specification of a
segment $S$ is described by the pair $(N,W)$, where $N$ is the
radical null-conditions ideal, and $W$ is a set of irreducible
(prime) polynomials on $K[\oa]$ representing non-null conditions
such that no prime component $N_i$ of the prime decomposition of
$N$ does contain any of the polynomials in $W$. We have $S=\V(N)
\setminus \V(h)$ with $h=\prod_{w\in W} w$. A red-specification
determined by $(N,W)$ is easily transformed into a {\em
diff-specification} $(N,M)$ with $N \subset M$ where $S=\V(N)
\setminus \V(M)$,  by considering the polynomial $h=\prod_{w\in W}
w$ and taking $M=\langle h \rangle + N$.

Let us denote ${\rm CGS}_1$ the output of BUILDTREE that consists
of a list of segments each represented by the three objects
$(B_i,N_i,W_i)$. Remember that each of these segments have
characteristic set of $\lpp$ of their bases $B_i$ that are
preserved by specialization on $S_i$. We say that a basis $G$ {\em
specializes well} to $(B,N,W)$, with $\lpp(G)=\lpp(B)$, if the
polynomials of $\overline{G}^{N}$ are proportional to the
polynomials of $B$, i.e. for each $g\in G$ there exist $f\in B$
and $\alpha,\beta \in W^*$ such that $\alpha \, \overline{g}^N =
\beta \, f$, where $W^* = \{k\prod_{i=1}^s w_i^{\lambda_i}: k\in
K, \lambda_i\in \mathbb{Z}_{\geq0},w_i\in W \}$.

Step ii), described in~\cite{Mo07}, selects the segments of ${\rm
CGS}_1$ with the same $\lpp$ that admit a common reduced basis
specializing well to the reduced Gr\"obner basis for every
specialization in the grouped segments. If Conjecture 7
in~\cite{Mo07} is true, the grouped segments form an intrinsic
partition of the parameter space. To perform that task the
algorithms DECIDE and GENIMAGE are used. The first one tests
whether one from two segments with the same $\lpp$ has already a
generic basis specializing to the other (this is the most frequent
case) or a sheaf exists and is necessary or whether possibly a
more generic basis must be found (by GENIMAGE). Whenever no
pre-image nor sheaf is found then both segments are not equivalent
and cannot be summarized. It can happen that instead of simple
polynomials the basis $B_i$ contains also sheaves of polynomials.
A sheaf $\{g_1,\dots,g_k\}$ is accepted in a basis of a segment
instead of a simple polynomial, whenever all the polynomials in
the sheaf specialize to the corresponding polynomial of the
reduced Gr\"obner basis of the specialized ideal or to $0$, and
some of the polynomials in the sheaf specialize to non-zero for
every $\oaa \in S_i$. As was shown in~\cite{Wi06}, it is necessary
to use sheaves for some over-determined systems if we want to
group all the segments admitting a common basis with the same
$\lpp$. We must notice that DECIDE algorithm also depends on the
semi-algorithm GENIMAGE to determine a polynomial $\tilde{f}$ that
specializes well to $f_1$ over $(N_1,W_1)$ and to $f_2$ over
$(N_2,W_2)$. Thus the canonicity of the results of the computation
of a MCCGS relies on GENIMAGE and the truthfulness of the
mentioned conjecture.

Let us denote the output of the second step ${\rm CGS}_2$. It will
be described by segments with a common basis $B_i$ and a set of
red-specifications:
\begin{equation}\label{grouping}
(B_i,\{(N_{i1},W_{i1}),\dots,(N_{ij_i},W_{ij_i})\}).\end{equation}
$S_i$ will now be the union of the segments determined by the
red-specifications $(N_{ik},W_{ik})$ for $k$ from 1 to $j_i$ .

Step iii) will be described in next sections. Its objective is to
give a canonical description of the union of the grouped segments
of step ii). In~\cite{Mo07} it was shown how a diff-specification
can be transformed into a {\em can-specification}. Here we will
prove that the union of red-specifications or their corresponding
diff-specifications can be transformed into a generalized
can-specification using what we call a $P$-tree. The idea is based
on Theorem 12 in~\cite{Mo07}. Let us give here a slightly
different formulation of it, more appropriate for the current
purposes.
\begin{theorem}\label{thm10}
\begin{enumerate}
  \item[]
  \item[{\rm i)}] Every diff-specification $S=\V(N)\setminus \V(M)$ admits a unique
  can-specification 
\begin{equation}\label{canspecform} S=\V(N) \setminus
\V(M)=\bigcup_{i} \left(\V(N_i) \setminus \left(\cup_{j}
\V(M_{ij}\right) \right),
\end{equation}
where ${\mathcal N}=\cap_i N_i$ and ${\mathcal M}_i=\cap_j M_{ij}$
are the irredundant prime decompositions over $A$ of the radical
ideals ${\mathcal N}$ and ${\mathcal M}_i$ respectively, where $
N_i \subsetneq M_{ij}$.
  \item[{\rm ii)}]  The Zariski closure over $\overline{K}^m$ verifies
    \[ \overline{S}=\overline{\bigcup_{i} \left(\V(N_i) \setminus \left(\cup_{j}
\V(M_{ij}\right) \right)} = \bigcup_i \V(N_i)=\V({\mathcal N}).\]
  \item[{\rm iii)}]  The can-specification verifies \[\V(N_i) \setminus \left(\cup_{j}
\V(M_{ij})\right)=S \cap \V(N_i).\]
 \item[{\rm iv)}] Given a diff-specification of $S$ the algorithm
 {\rm DIFFTOCANSPEC}~{\rm \cite{Mo07}} builds its can-specification.

\end{enumerate}
\end{theorem}
The need of having a canonical description of the intrinsic
segments comes from the need of comparing different outputs for
the same problem, and also to have a final simple description of
the segments.

\section{Adding segments}\label{addingsection}
We tackle now the third step of MCCGS, i.e. the description of the
union of the segments in a canonical form. We start with segments
of the form~(\ref{grouping}). The red-specifications $(N,W)$ can
be transformed into diff-specifications $(N,M)$, as explained in
Section~\ref{SecPreliminaries}, so we are attained with the
obtention of a canonical representation for the addition of
diff-specifications. We cannot assume that the simple form given
by formula~(\ref{canspecform}) will be sufficient. A more complex
constructible set will be formed grouping all the segments
$S_{ik}$ for $1 \le k \le j_i$.

Thus we generalize the concept of canonical specification given in
\cite{Mo07}:
\begin{definition}[P-tree] \label{defPtree}
A P-tree is a rooted directed tree such that
\begin{enumerate}\renewcommand*{\labelenumi}{\roman{enumi})}
  \item  the nodes are prime ideals over $A$ except the root, denoted $r$,
  \item  when $P \rightarrow Q$ is an arc then $P \subsetneq
  Q$,
  \item  the children of a node are a set of irredundant prime
  ideals over A, (whose intersection form a radical ideal).
\end{enumerate}
By definition the root level is $0$.
\end{definition}
\begin{definition}[C-tree]\label{defCtree}
To any P-tree we associate an isomorphic C-tree by changing every
node $P$ to a subset of $\overline{K}^m$ denoted $C(P)$ by the
following recursive procedure:
\begin{enumerate}\renewcommand*{\labelenumi}{\roman{enumi})}
  \item  if $P$ is a leaf (terminal vertex) then
  $C(P)=\V(P)$,
  \item  if $P$ is an inner node different from the root and $P_1,\dots,P_d$ are
  its children, then
  \begin{equation} \label{gencanspecform} C(P)=\V(P) \setminus \left( C(P_1) \cup \dots \cup C(P_d)
  \right)\end{equation}
  \item  if $P_1,\dots,P_d$ are
  the children of the root vertex $r$ then
  \[ C(r) = C(P_1) \cup \dots \cup
  C(P_d). \]
\end{enumerate}
\end{definition}
Note that for $C(r)$ the parity of the vertex-level acts
additively for odd level vertices and as a subtraction for even
level vertices. (See example~\ref{exampletree} below).
\begin{definition}[Generalized canonical specification]\label{defGCS}
A generalized canonical specification (GCS) of a set $S$ is a
P-tree such that $S=C(r)$ satisfying, for every node $P$ at level
$j$, the following condition:
\begin{equation} \label{gencanspecprop} C(P)=\V(P) \cap B\end{equation}
 where $B=S$ for $j$ odd and $B=\overline{K}^m\setminus S$ for $j$ even.
\end{definition}
\begin{example}\label{exampletree}
To clarify the definition suppose that we want to describe the set
$S_1$ of the $\Real^3$-space with coordinates $a,b,c$ consisting
of the planes $a=0$ and $b=-1$ except the lines $a=b=0$ and
$a=c=0$ plus the point $O(0,0,0)$. We can express $S_1$ as
\[ S_1= \left( \left( \V(a) \cup \V(b+1) \right)
\setminus \left( \V(a,b) \cup \V(a,c) \right) \right) \cup
\V(a,b,c)\]
But there exist many other possible determinations of
this set. If we want to obtain the GCS of $S_1$ we must write
$S_1$ in the form
\[ S_1 = \left(\V(a) \setminus
\left(\left( \V(a,b) \setminus  \V(a,b,c) \right) \cup \left(
\V(a,c) \setminus \left( \V(a,b,c) \cup \V(a,b+1,c) \right)\right)
\right)\right) \cup \V(b+1)
\]
\begin{figure}
\begin{center}
\includegraphics[width=13cm]{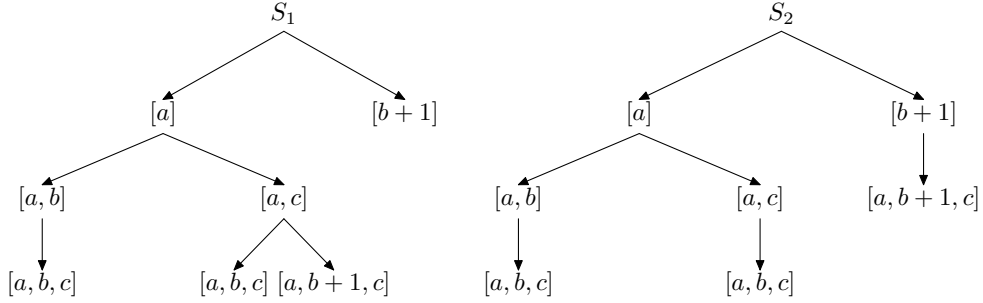}
\caption{\label{examplegcstree} Trees representing the sets $S_1$
and $S_2$ in generalized can-specification.}
\end{center}
\end{figure}
This formula can be represented by the tree associated to $S_1$
shown in figure~\ref{examplegcstree}. Notice that we must include
the point $\V(a,b+1,c)$ under the branch of $\V(a,c)$, as this
point belongs to $S_1$ and condition (\ref{gencanspecprop})
requires it to belong also to $C(a)$. The interest of that
representation lies in the fact that it is unique as we prove in
Theorem \ref{gencanspecthm} below.

Consider now the set $S_2=S_1 \setminus \V(a,b+1,c)$. In order to
preserve property (\ref{gencanspecprop}) of the GCS definition,
the P-tree associated to $S_2$ will be modified from the P-tree
associated to $S_1$ by eliminating the point under the variety
$\V(a,c)$ and setting it under the variety $V(b+1)$. The new tree
is also shown in Figure \ref{examplegcstree}. These examples
should clarify the definition of GCS to obtain canonicity of the
description by preserving condition (\ref{gencanspecprop}).
\end{example}
\begin{theorem}\label{gencanspecthm}
A subset $S\subset \overline{K}^m$ defined by a {\rm GCS} has the
following properties:
\begin{enumerate}
  \item[]
  \item[{\rm i)}]  For every vertex $P$, except for the root,
  \[ \overline{ C(P)} = \V(P)\]
 where, as usual, the Zariski closure is taken over $\overline{K}^m$.
  \item[{\rm ii)}]  For the root vertex $r$
\[\overline{S}=\overline{C(r)}=\bigcup_{i=1}^d \V(P_{i})\]
where the $P_i$'s are the children vertices of $r$.
  \item[{\rm iii)}]  $S$ has a unique {\rm GCS} decomposition.
\end{enumerate}
\end{theorem}
\begin{proof}
\begin{enumerate}\renewcommand*{\labelenumi}{\roman{enumi})}
   \item The inclusion $\subseteq$ is obvious as $C(P) \subseteq \V(P)$. To prove the equality we have
 \[C(P)=\V(P) \setminus \bigcup_{i=1}^d C(P_i)
   \supseteq \V(P)
 \setminus  \bigcup_{i=1}^d \V(P_i).\]
  Consider the closure of the above formula and
  apply Theorem~\ref{thm10}~(ii). The result follows.
   \item  Is an immediate consequence of i).
  \item  To prove the uniqueness we proceed by induction on $d$. For
  $d=1$, the tree is formed by the root $r$ and a set of children nodes
  forming an
  irredundant prime decomposition of the radical ideal defining $S$, by
  Definition~\ref{defPtree}~iii).
  Thus, in this segment the P-tree is unique.

 Assume now by induction hypothesis the
 uniqueness of the GCS for every P-tree of maximum depth less
 than $d$ and let us prove, as a consequence, the uniqueness also for
 depth $d$. Let $S$ be defined by a P-tree of maximal depth $d$ representing a
 GCS.
 By part (ii) of the Theorem we have
 \[ \overline{S}=\bigcup_{i} \V(P_{i}) =
 \V(\displaystyle\cap_{i} P_{i}),
\]
 where the $P_{i}$'s form the unique irredundant prime
 decomposition over $A$ of the radical ideal $\cap_{i} P_{i}$ defining $\overline{S}$ by Definition~\ref{defPtree}~(iii).
 Thus they are uniquely determined. Denoting $P_{ij}$ the children of $P_i$, by
(\ref{gencanspecprop}), we have
 \begin{equation}\label{GCSPi}
  C(P_i)=\V(P_i) \setminus \displaystyle \bigcup_{j=1}^{d_i} C(P_{ij}) = \V(P_i) \cap S
 \end{equation}
 showing that $C(P_i)$ is also uniquely determined.
 Set $S_i$ for the subtracting set
 \begin{equation}\label{defCi} S_i=\displaystyle \bigcup_{j=1}^{d_i} C(P_{ij}).\end{equation}
 As $S_i \subseteq \V(P_i)$,
 $S_i$ is also uniquely defined by (\ref{GCSPi}). By Definition \ref{defGCS}, formula (\ref{gencanspecprop}), we have
\[
C(P_{ij}) = \V(P_{ij}) \cap \left( \overline{K}^m \setminus S \right),\\
 \]
Thus
\[S_i=\displaystyle \bigcup_{j=1}^{d_i} C(P_{ij})=\left( \displaystyle \bigcup_{j=1}^{d_i} \V(P_{ij}) \right) \cap
\left( \overline{K}^m \setminus S \right) \] and so
\begin{equation}\label{propCij} C(P_{ij})=\V(P_{ij}) \cap S_i\end{equation}
By the ascending chain condition for the ideals in the branches
and condition (\ref{gencanspecprop}) for the P-tree of $S$,
equation (\ref{propCij}) ensures that condition
(\ref{gencanspecprop}) is also respected for the subtree of $S_i$,
whose root vertex is given by (\ref{defCi}). Thus the subtree of
$S_i$ also forms a GCS of $S_i$ with depth less than $d$. By the
induction hypothesis it is uniquely determined and so does the
complete P-tree of $S$.
\end{enumerate}
\end{proof}

\section{The MCCGS algorithm} \label{algorithmsection}
Given an ideal $I$ and the monomial orders $\succ_{\ox}$ for the
variables and $\succ_{\oa}$ for the parameters, the following
sequence of algorithms build up the $P$-tree $T$ corresponding to
the Minimal Canonical Comprehensive Gr\"obner System associated to
$I$ and $\succ_{\ox}$. We describe them in descendent design.
\begin{algorithm}[ht]
 tree $T$ $\leftarrow$ {\bf MCCGS}$(B,\succ_{\ox},\succ_{\oa})$ \\
 {\it Input:} $B$ a basis of the parametric polynomial ideal $I$ and monomial orders $\succ_{\ox},\succ_{\oa}$.\\
 {\it Output:} $T$ a tree containing the minimal canonical comprehensive Gr\"obner system associated to
 $I$.\\
 \begin{algorithmic}[]
   \STATE {$T_0:=${\bf BUILDTREE}$(B,\succ_{\ox}, \succ_{\oa})$}
   \STATE {$S:=${\bf SELECTCASES}$(T_0)$}
   \STATE {$T:=${\bf GENCANTREE}$(S)$}
 \end{algorithmic}
\end{algorithm}
\begin{algorithm}[ht]
 set of pairs $S$ $\leftarrow$ {\bf SELECTCASES}$(T_0)$ \\
 {\it Input:} $T_0$ a BUILDTREE discussion tree whose terminal vertices shape a CGS with red-specifications. \\
 {\it Output:} $S$ a finite set of pairs of the form $(B_i,\{(N_{i1},W_{i1}),\dots,(N_{ij_i},W_{ij_i})\})$ taken from the CGS associated to $T_0$. \\
 \begin{algorithmic}[]
   \STATE $G:=\{(B_1,N_1,W_1),\dots,(B_r,N_r,W_r)\}$ \COMMENT{the CGS associated to $T_0$}
   \STATE {$S:=\emptyset$}
   \WHILE {$G\not=\emptyset$}
     \STATE {Let $(B,N,W)$ be the first element of $G$}
     \STATE {$B_0:=B;\  N_0:=N;\  W_0:=W;$}
     \STATE {$l:=\{(N_0,W_0)\}$}
     \STATE {$G:=G\setminus \{(B_0,N_0,W_0)\}$}
     \FORALL {$(B',N',W')\in G$ such that $\lpp(B)=\lpp(B')$}
       \STATE {$p:=0$}
       \FORALL {$f \in B$ {\bf while} $p \neq \hbox{\bf false}$}
         \STATE {Let $f'\in B'$ be such that $\lpp(f)=\lpp(f')$}
         \STATE {$p:=\hbox{\bf DECIDE}(f,N,W,f',N',W')$}
         \IF {$p \neq \hbox{\bf false}$}
           \STATE {Substitute $f$ by $p$ in $B_0$}
         \ENDIF
       \ENDFOR
       \IF {$p \neq \hbox{\bf false}$}
         \STATE {$l:=l\cup \{(N',W')\}$}
         \STATE {$B:=B_0;\  N:=N\cap N';\  W:= W\cap W';$}
       \ENDIF
     \ENDFOR
     \STATE {$S:=S\cup \{(B,l)\}$}
     \STATE {$G:=G\setminus \{(B',N',W')\in G \text{ such that } (N',W') \in l\}$}
  \ENDWHILE
 \end{algorithmic}
\end{algorithm}

MCCGS uses BUILDTREE (see \cite{Mo02, MaMo06}) to build up the
discussion tree $T_0$ containing a CGS whose segments are
expressed as red-spe\-ci\-fi\-ca\-tions. Then SELECTCASES takes
$T_0$ as input and classifies the segments from the CGS associated
to $T_0$ into pairs of the form $(B_i,l_i)$, where $l_i$ is a set
of red-specifications $\{(N_{i1},W_{i1}),\dots (N_{ij_i},W_{ij_i})
\}$ whose corresponding bases have been generalized by the same
basis $B_i$. Afterwards, MCCGS calls the new algorithm GENCANTREE
to finally obtain the MCCGS associated to the initial ideal and
term order.

\begin{algorithm}[ht]
 tree $T$ $\leftarrow$ {\bf GENCANTREE}$(S)$ \\
 {\it Input:} $S$ a finite set of pairs of the form $(B_i,\{(N_{i1},W_{i1}),\dots,(N_{ij_i},W_{ij_i})\})$. \\
 {\it Output:} the canonical tree $T$ associated to $S$.\\
 \begin{algorithmic}[]
   \STATE {initialize $T$}
   \FOR{$1\leq i\leq \sharp S $}
     \STATE {Create $u_i$ a new vertex in $T$ hanging from root}
     \STATE {store $B_i$ in $u_i$}
     \STATE $l:=\{(N_{i1},W_{i1}),\dots,(N_{ij_i},W_{ij_i})\}$
     \COMMENT{red-specifications associated to $B_i$}
     \STATE {$\overline{T}:=${\bf GCS}$(l)$}
     \STATE {hang $\overline{T}$ from $u_i$}
   \ENDFOR
 \end{algorithmic}
\end{algorithm}
\begin{algorithm}[ht]
 tree $\overline{T}$ $\leftarrow$ {\bf GCS}$(l)$ \\
 {\it Input:} $l$ a finite set of red-specifications\\
 {\it Output:} a tree containing the Generalized Can-Specification associated to the addition of segments in
 $l$.\\
 \begin{algorithmic}[]
   \STATE {initialize tree $\overline{T}$ with root $r$}
   \STATE {set $P_r:=\phi$}
   \FORALL {pairs $(N,W)\in l$}
     \STATE {$\overline{T}:=${\bf ADDCASE}$((N,W),r,\overline{T})$}
   \ENDFOR
 \end{algorithmic}
\end{algorithm}

GENCANTREE uses GCS algorithm to build the $P$-tree corresponding
to the generalized canonical specification of the addition of
segments. GCS algorithm begins by setting the ideal $\{0\}$ at the
root of new tree $\overline{T}$ and calls iteratively the
recursive algorithm ADDCASE. It must be noted that there are two
kinds of nodes, namely odd level vertices and even level vertices,
that are treated differently by ADDCASE. ADDCASE uses two
auxiliary algorithms: DIFFTOCANTREE (a minor transformation of
DIFFTOCANSPEC) converts a diff-specification into a $P$-tree
containing the associated can-specification, and SIMPLIFYSONS just
makes the suitable simplifications.

\begin{figure}
\begin{center}
\includegraphics{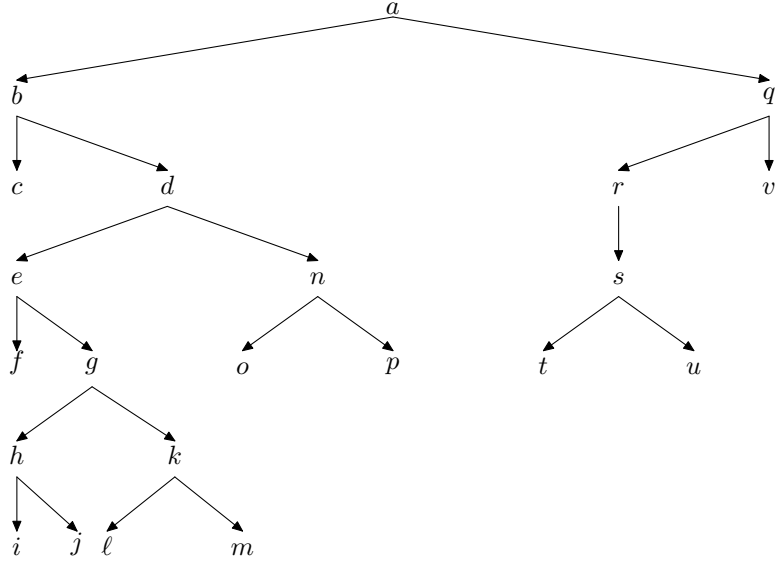}
\caption{\label{addsegmentgraf} The action of ADDCASE.}
\end{center}
\end{figure}

At the first iteration ADDCASE stores under root the $P$-tree of
the unique canonical specification associated to
$(N_{i1},W_{i1})$. Then, to add each further red-specification
$(N_{ik},W_{ik})$, ADDCASE executes itself recurrently in
post-order at the even level vertices $u \in \overline{T}$ and
adds the can-specification associated to $(N_{ik},W_{ik})$
contained in $\V(P_u)$. For example, in
figure~\ref{addsegmentgraf} it would act successively on the
vertices
\[ c,\ f,\ i,\ j,\ \ell,\ m,\ g,\ o,\ p,\ d,\ t,\ u,\ r,\ v,\ a.\]

\begin{algorithm}[ht]
 (bool, tree $\overline{T}$) $\leftarrow$ {\bf ADDCASE}$((N,W),u,\overline{T})$ \\
 {\it Input:} $(N,W)$ a red-specification , $u$ the current vertex in $P$-tree $\overline{T}$.\\
 {\it Output:} $false$ if $(N,W)$ is not to be added to parent vertices, $true$
 otherwise. It also returns current tree $\overline{T}$.\\
 \begin{algorithmic}
   \IF {$u$ is not terminal}
     \STATE {$test:=\hbox{\bf treu}$}
     \FORALL {$v \in \hbox{\rm children}(u)$}      
       \FORALL {$w \in \hbox{\rm children}(v)$}   
         \IF {{\bf ADDCASE}$((N,W),w)=\hbox{\bf false}$}
           \STATE{$test:=\hbox{\bf false}$ }
         \ENDIF
       \ENDFOR
       \STATE {$\overline{T}:=\hbox{\bf SIMPLIFYSONS}(v,\overline{T})$}
     \ENDFOR
   \ELSE
     \STATE {$test:=\hbox{\bf true}$}
   \ENDIF
   \IF {$test=\hbox{\bf true}$}
     \STATE {$h:=\prod_{w\in W} w$}
     \STATE $(R,S):=(N+P_u,N+\langle h\rangle + P_u)$ \COMMENT{diff-specification associated to $(N,W)$ in $\V(P_u)$ }
     \STATE {$t:=\hbox{\bf DIFFTOCANTREE}(R,S)$}
     \STATE {hang $t$ from $u$}
     \IF {parent$(u)$ exists {\bf and} $P_{\text{parent}(u)} \subseteq R $}
     \STATE {$test:=\hbox{\bf false}$}
     \ENDIF
   \ENDIF
 \end{algorithmic}
\end{algorithm}
\begin{algorithm}[ht]
 tree $\overline{T}$ $\leftarrow$ {\bf SIMPLIFYSONS}$(v,\overline{T})$ \\
 {\it Input:} $v$ a vertex at odd level of tree $\overline{T}$ where to start the simplifications. \\
 {\it Input:} The tree after simplifications\\

 \noindent {\it Description}:\\
 SIMPLIFYSONS just simplifies the subtree under $v$ on the global $\overline{T}$ in order to not having
 cancellations nor inclusions between the children of $v$. Let $P$ be the prime stored in vertex $v$. The simplification
 is performed as follows: \\
 Check that there is no $P_i$ child of $P$ such that $P_i=P_{ij}$. And if any, hang to $P$ all subtrees
 descendant from $P_{ij}$ and drop both $P_i$ and $P_{ij}$ from $\overline{T}$.

 Then check whether there is any pair of children of $P$,
 $P\rightarrow P_{i}$, $P\rightarrow P_{j}$, such that $P_{i} \subseteq P_{j}$. If so, drop subtree hanging
 from $P_{j}$ and also vertex $P_j$.
 \begin{algorithmic}[]
 \FORALL {$v\in$ children$(u)$}
   \IF {$P_v=P_{\text{child}(v)}$}
     \STATE{hang from $u$ all subtrees under child$(v)$}
     \STATE{drop $v$ and child$(v)$ from $\overline{T}$}
   \ENDIF
 \ENDFOR
 \IF {there $\exists \  v,w \in $children$(u)$ such that $P_v\subseteq P_w$}
   \STATE {drop subtree with root $w$ from $\overline{T}$}
 \ENDIF
 \end{algorithmic}
\end{algorithm}
\begin{algorithm}[]
 tree $t$ $\leftarrow$ {\bf DIFFTOCANTREE}$((I,J))$ \\
 {\it Input:} $(I,J)$ a diff-specification. \\
 {\it Output:} a tree structure containing the Can-Specification of $\mathbb{V}(I) \setminus
 \mathbb{V}(J)$.\\
 \begin{algorithmic}[]
   \STATE {initialize local tree $t$}
   \STATE {$\{P_i\}:=${\bf PRIMEDECOMP}$(I)$}
   \FORALL {$P_i$}
     \IF {$P_i \not= \sqrt{J+P_i}$}
       \STATE {store the $P_i$ as the children of root in $t$}
       \STATE {$\{P_{ij}\}:=${\bf PRIMEDECOMP}$(J+P_i)$}
       \STATE {store the $P_{ij}$ as the children of $P_i$ in $t$}
     \ENDIF
   \ENDFOR
 \end{algorithmic}
\end{algorithm}

Thus, before acting on an even vertex $u \in \overline{T}$, the
algorithm must have acted on all its even descendants. Therefore,
if an even level descendant $v$ verifies that $N_{ik}\supseteq
P_v$, then the can-specification associated to $(N_{ik},W_{ik})$
must have been completely hung under $v$. In this case the
\emph{test} variable will contain \emph{false} and thus
DIFFTOCANTREE for current $(N_{ik},W_{ik})$ will not act on $P_u$
nor on any of its ascendant vertices. We must also remind that the
ideals associated to the paths in $\overline{T}$ starting from
root form ascending chains of prime ideals. Thus, whenever
\emph{test} is \emph{false}, the condition cited above will also
hold for all vertices placed between $u$ and $v$, even the odd
level ones, i.e. for all $w\in \overline{T}$ descendent of $u$ and
ascendant of $v$, $N_{ik} \supseteq P_w$.

This way, ADDCASE completes current P-tree $\overline{T}$ to a new
tree such that for every odd level vertex $u$ with prime ideal
$P_u$, all points in $\V(P_u) \cap \left( \V(N_{ik})\setminus
\V(h_{ik})\right)$ (where $h_{ik}=\prod_{w\in W_{ik}}w$) are in
$C(P_u)$, as required.

Nevertheless in the new tree completed by ADDCASE it could happen
that $P_u+N_{ik}=P_u$ for some even level vertex $u$, which would
cause that $P_u$ and its unique child $P_{\text{child}(u)}$
coincide. If so, SIMPLIFYSONS takes the subtree under child$(u)$,
slips it upwards hanging it from parent$(u)$ and eliminates both
vertices $u$ and child$(u)$ from the tree. When this action is
performed, it could also happen that some set of current even
level siblings do not preserve the prime decomposition
irredundancy property, as some lifted primes can contain some of
their sibling vertices, i.e. $\exists v_1,v_2 \in
\text{children}(u)$ such that $P_{v_1} \subseteq P_{v_2}$ for $u$
an even level vertex in $\overline{T}$. SIMPLIFYSONS algorithm
also detects these cases and eliminates the subtrees hanging from
$v_2$ as well as $v_2$. Though, the action of SIMPLIFYSONS will
restore the GCS-condition property of the tree.

Note: For algorithmic reasons, all paths starting from root vertex
in a $P$-tree will be of even length. Thus for odd length
branches, the algorithm will add a new vertex $[1]$ at the end.

The above described algorithms build the complete MCCGS of the
initial ideal. The following theorem states that GCS algorithm
builds the generalized can-specification (GCS) associated to the
set of the corresponding diff-specifications:

\begin{theorem}
Given a finite list of pairs $l=\{(N_{ik},W_{ik}): k=1,\dots, M\}$
of red-specifications, {\rm GCS}$(l)$ computes the $P$-tree
associated to the generalized can-specification determining the
constructible set
$$\bigcup_{k=1}^M \V(N_{ik})\setminus \V(\prod_{w\in W_{ik}}w).$$
\end{theorem}
\begin{proof}
Let $S=\bigcup_{k=1}^M \V(N_{ik})\setminus \V(\prod_{w\in
W_{ik}}w)$. The proof is done by induction on $M$, the number of
red-specifications to be added.

For $M=1$, GCS uses DIFFTOCANTREE just once and, by Theorem 1
(iv), it builds up the unique can-specification in tree
$\overline{T}$. Thus $\overline{T}$ is a $P$-tree such that
$C(\overline{T})=\V(N_{i1})\setminus\V(\prod_{w\in W_{i1}}w)$.

By induction hypothesis, assume now that after the $M-1$ iteration
of ADDCASE the GCS tree of the $M-1$ red-specifications has been
built and let $\tilde{T}$ be this tree, which is a $P$-tree such
that $C(\tilde{T})=\bigcup_{k=1}^{M-1} \V(N_{ik})\setminus
\V(\prod_{w\in W_{ik}}w)$ and such that every vertex $u \in
\tilde{T}$ holds that $C(P_u)=\V(P_u)\cap C(\tilde{T})$. We shall
prove that the $M$-th iteration will build the GCS tree of $S$.

Let us describe how the recursive ADDCASE algorithm acts on
$\tilde{T}$ adding $\V(N_{iM})\setminus\V(\prod_{w\in W_{iM}}w)$.
Denote by $\Lambda(u)$ the operation on an even level vertex $u$
that hangs to it the tree associated to the can-specification of
$(N_{iM},W_{iM})$ contained in $\V(P_u)$ (i.e.
$\V(N_{iM}+P_u)\setminus \V(N_{iM}+P_u+\langle \prod_{w\in
W_{iM}}w \rangle)$ whenever it can be hung and returns $false$ or
$true$ depending on whether parent$(P) \subseteq N_{iM}$ or not,
respectively. So it hangs the points $\V(P_u) \cap (
\V(N_{iM})\setminus\V(N_{iM}+\langle \prod_{w\in W_{iM}}w
\rangle))$, and thus $C(P_u)=\V(P_u) \cap S$.

$\Lambda(u)$ is applied recursively in post-order. If $\Lambda(u)$
returns $false$ at some even level vertex $u$, the whole set
$\V(N_{iM})\setminus\V(N_{iM}+\langle \prod_{w\in W_{iM}}w
\rangle)$ has been hung under $u$ and thus, as $u$ is even,
$C(\text{father}(u))\supset \V(N_{iM})\setminus\V(N_{iM}+\langle
\prod_{w\in W_{iM}}w \rangle)$. Then $\Lambda$ will not be applied
to any of its ascendant vertices because $C(\overline{T})=
C(\tilde{T})\cup \left( \V(N_{iM})\setminus\V(N_{iM}+\langle
\prod_{w\in W_{iM}}w \rangle)\right)$.

If $\Lambda(u)$ returns $true$ for all $u \in \tilde{T}$, which
means that $\V(N_{iM})\setminus\V(N_{iM}+\langle \prod_{w\in
W_{iM}}w \rangle)$ has not completely been hung under root, then
the $P$-tree corresponding to the red-specification
$(N_{iM},W_{iM})$ computed by DIFFTOCANTREE will be hung from
root. Thus, we finally have that $C(\overline{T})=
C(\tilde{T})\cup \left( \V(N_{iM})\setminus\V(N_{iM}+\langle
\prod_{w\in W_{iM}}w \rangle)\right)$.

This way, GCS algorithm obtains, as SIMPLIFYSONS ensures, a
$P$-tree $\overline{T}$ such that for every node $v \in
\overline{T}$ holds that $C(P_v)=\V(P_v)\cap C(\overline{T})$ and
$C(\overline{T})=S$.
\end{proof}

Furthermore, GENCANTREE algorithm performs a GCS computation for
each list of segments whose associated reduced Gr\"obner bases
specialize properly, obtaining a tree for which the subtrees
hanging from the root correspond to the generalized
can-specifications of the lists configuring a partition of the
parameter space. Thus, MCCGS algorithm performs the discussion and
obtains the Minimal Canonical Comprehensive Gr\"obner System
stored in the output tree $T$.

\begin{example}\label{conicex}[Singular points of a conic]
The general equation of a conic can be reduced by a suitable
change of variables to the form
\[f\equiv x^2+b y^2+2 c x y+d x=0.\] To study its singular points
consider the system of equations
\begin{figure}
\begin{center}
\includegraphics{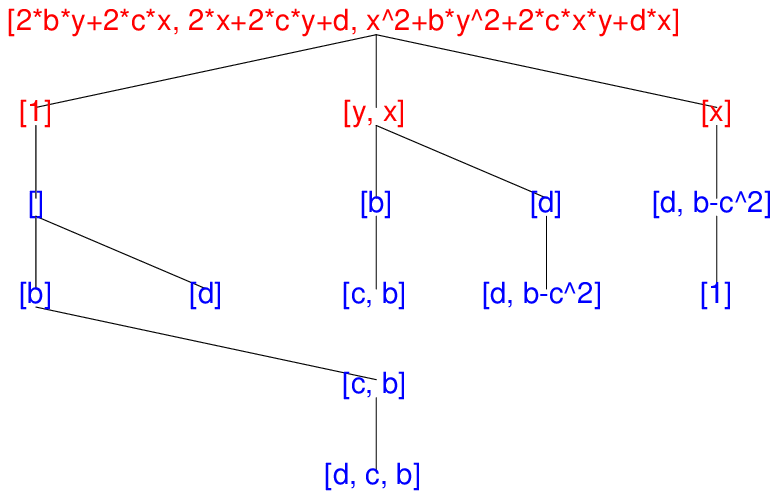}
\caption{\label{conicsingpointtree} MCCGS for the singular points
of a conic system.}
\end{center}
\end{figure}
\[S:=\left[f,\frac{\partial{f}}{\partial{x}},\frac{\partial{f}}{\partial{y}}\right], \]
and apply MCCGS algorithm to $S$ using $\lex(x,y)$ and
$\lex(b,c,d)$ for variables and parameters respectively. The
result is shown in Figure~\ref{conicsingpointtree}. The
interpretation of the output tree is the following.

There are three different segments: The generic case with $\lpp$
set $[1]$ where the conic has no singular points, the segment with
$\lpp$ set $[y,x]$ corresponding to a single singular point in the
conic, and the segment with $\lpp$ set $[x]$ corresponding to a
solution with one degree of freedom, where the conic is a double
line. The conditions over the parameters given by the trees are to
be interpreted in the following way:\vspace{3mm}

\begin{center}
\begin{tabular}{||c|l|l||}
 \hline \hline
 $\lpp$ & Basis & Description \\
 \hline
 $[1]$ & $[1]$ & $\Compl^3 \setminus \left( \left( \V(b)
\setminus (\V(c,b) \setminus \V(d,c,b)) \right) \cup \V(d) \right)
 $ \\
 \hline
 $[y,x]$ & $[2 c y+d, x]$ & $\left(\V(b) \setminus \V(c,b) \right) \cup  \left( \V(d) \setminus
 \V(d,b-c^2)\right)
  $ \\
  \hline
  $[x]$ &  $[x+cy]$ & $ \V(d,b-c^2) $ \\
 \hline\hline
\end{tabular}
\end{center}\vspace{3mm}

Figure~\ref{conicsingpointvars} shows the geometrical description
of the partition of the parameter space provided by the three
segments. The generic segment occurs in the whole 3-dimensional
space except the two planes $\V(c)$ and $\V(d)$ plus the line
$\V(c,b)$ except the point $(0,0,0)$. The one-singular point
segment occurs in the two planes $\V(c)$ and $\V(d)$ except both
the line $\V(c,b)$ and the parabola $\V(d,b-c^2)$. Finally the
double line occurs on the parabola $\V(d,b-c^2)$.
\begin{figure}
\begin{center}
\includegraphics{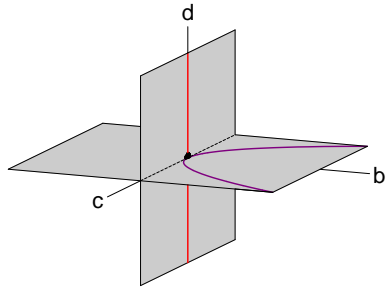}
\caption{\label{conicsingpointvars} Geometrical description of the
MCCGS for the singular points of a conic system.}
\end{center}
\end{figure}

\end{example}

\section{Applications} \label{applicationsection}
\begin{figure}
\begin{center}
\includegraphics{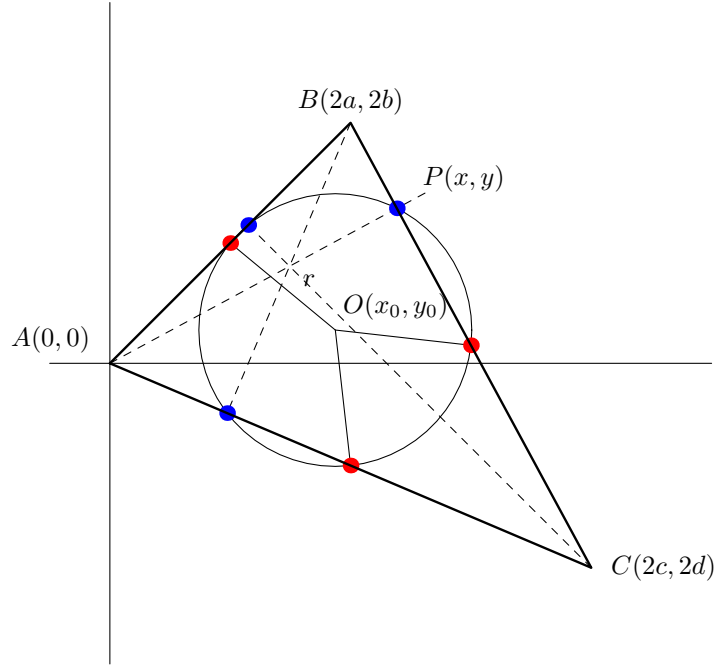}
\caption{\label{ninepointcircleth} Nine points circle Theorem.}
\end{center}
\end{figure}

\begin{figure}
\begin{center}
\includegraphics{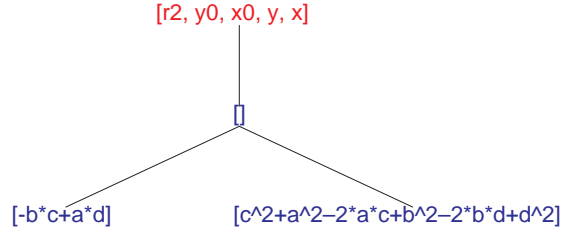}
\caption{\label{Figninepointcircletree} Generic case for $HT$ in
the nine points circle Theorem.}
\end{center}
\end{figure}

We use now the algorithm to prove part of the 9 points circle
Theorem on a triangle. It states: {\em For every triangle, the
circle through the three middle points of the sides is also
incident with the height feet}. To prove it, and also to obtain
supplementary hypotheses if needed, consider a triangle with
vertices at the points $A(0,0)$, $B(2a,2b)$ and $C(2c,2d)$ and
denote $P(x,y)$ the height foot from $A$ (see Figure
\ref{ninepointcircleth}).
The first set of hypotheses are the equations of the side $BC$ and
the height from $A$ defining the point $P(x,y)$:
\[\begin{array}{ll}
 h_1: & (b-d)x+(c-a)y+2ad-2bc=0\\
 h_2: & (c-a)x+(b-d)y=0
\end{array}
\]
Denote $r$ and $(x_0,y_0)$ the radius and the center of the circle
through the three middle points $(a,b)$, $(c,d)$ and $(a+c,b+d)$.
Its equation will be $(x-x_0)^2+(y-y_0)^2-r^2=0$. So we have the
three new hypotheses:
\[\begin{array}{ll}
 h_3: & (a-x_0)^2+(b-y_0)^2-r^2=0\\
 h_4: & (c-x_0)^2+(d-y_0)^2-r^2=0\\
 h_5: & (a+c-x_0)^2+(b+d-y_0)^2-r^2=0
\end{array}
\]
The thesis of the theorem is that the circle is incident with the
point $P(x,y)$, thus that the polynomial
\[f=(x-x_0)^2+(y-y_0)^2-r^2\]
is zero as a consequence of the hypotheses. The first to do is
searching for the solutions of the system $HT = \langle
h_1,h_2,h_3,h_4,h_5,f \rangle $. Thus we call
 \[\hbox{\rm mccgs}(HT,\grevlex(x,y,x_0,y_0,r_2),\lex(a,b,c,d)),\]
where we set $r_2=r^2$. We obtain a canonical tree with nine
cases. But only two cases are really interesting. The first one is
the generic case (see Figure~\ref{Figninepointcircletree}) for
which the $\lpp$ are $[r_2,y_0,x_0,y,x]$ showing that for
parameter values not in $\V(ad-bc) \cup \V((a-c)^2+(b-d)^2)$ there
exists a unique solution. For the real case it is sufficient to
consider $ad-bc \ne 0$, as the real part of the second variety is
inside the first one. The second interesting case is the case with
basis $[1]$ where no solution exists. The corresponding tree shows
that it covers both varieties $\V(ad-bc) \cup \V((a-c)^2+(b-d)^2)$
except for very special cases corresponding to degenerate
triangles. Thus we have proved that the theorem is true whenever
$ad-cb\neq 0$.
We can also go further and ask if the thesis is a real consequence
of the hypotheses, i.e. if $f$ belongs to the radical of the
hypotheses ideal ${H} = \langle h_1,h_2,h_3,h_4,h_5 \rangle $
whenever $ad-bc \neq 0$. To test this we must have
\[ {HT}_1= \langle h_1,h_2,h_3,h_4,h_5,1-wf \rangle =\langle 1 \rangle\]
i.e. the Gr\"obner basis of ${HT}_1$ is $[1]$. We call now
 \[\begin{array}{ll}
 \hbox{\rm mccgs}&([h_1,h_2,h_3,h_4,h_5,1-wf],\grevlex(w,x,y,x_0,y_0,r),\\
 & \lex(a,b,c,d),notnull=\{ad-bc\}),
 \end{array}\]
and the result is a unique case with basis $[1]$. Thus effectively
$f$ belongs always to the ideal of the hypothesis whenever
$ad-bc\neq 0$.

\section{Comparison of algorithms}\label{compare}
The CGS of a parametric ideal $I$ can have very different
properties as commented in section~\ref{intro}. For example
\begin{enumerate}\renewcommand*{\labelenumi}{\roman{enumi})}
  \item  the subsets
$S_i$ of the parameter space $\overline{K}^m$ in which the CGS are
divided can be very different, they can contain different number
of segments, they can overlap, and so on;
  \item  a CGS can contain incompatible segments;
  \item  the basis $B_i$ can be reduced or not;
  \item  even when a given algorithm does not theoretically  ensure some
property it can however hold it experimentally in most examples.
\end{enumerate}
So it is quite difficult to make automatic comparisons of the
outputs.

There are three available known implemented methods for obtaining
a GCS:
\begin{enumerate}\renewcommand*{\labelenumi}{\roman{enumi})}
  \item Weispfenning CGB implemented by~\cite{DoSeSt06}
 in Reduce.
  \item Suzuki-Sato SACGB implemented in
  Risa/Asir\footnote{There exist also a preliminary Maple version but it is not yet fully developed.},
  \item Montes MCCGS implemented in Maple 8 by M.
  Manubens in the DPGB library 7.0.
\end{enumerate}

Even though we use some criteria to evaluate them:
  correctness of the results,
  existence of incompatibilities,
  existence of overlaps,
  number of segments,
  whether the $S_i$ form a partition,
  whether or not specializations preserve the $\lpp$'s of the bases,
  reduction of the bases,
  theoretical canonicity ensured,
  theoretical minimality ensured,
  execution time.

Although it is not possible to evaluate Weispfenning's CCGB
algorithm in practice because it has not been implemented, we can
analyze its theoretical features. The canonicity of CCGB comes
from the use of primary decompositions over the conditions, but
the method is not dichotomic and so the segments are not disjoint.
As its objective is to obtain a canonical CGB, the bases of the
corresponding CGS are faithful and therefore not reduced, so
specializations do not preserve their $\lpp$. Furthermore as the
segments are not disjoint, minimality does not hold.

For the comparisons with the implemented methods, we have used a
Pentium(R) D CPU 3.00 GHz, 1.00 GB RAM for the computations and
tested different examples using the above implementations.

We have not been able to obtain CGB Reduce in time for these
comparisons, so we could only test some very simple executions
with a demo version. To what we have experimentally observed, it
gives a partition of the parameter space containing quite more
segments than MCCGS. The bases are faithful, which is interesting
to compute a CGB, but do not give direct information on the type
of solutions, as these bases are not reduced. It seems to be very
efficient but the provided results are difficult to be
interpreted. In the future we will make a more precise analysis.

SACGB is a very simple and interesting algorithm based on
Kalkbrenner's theorem for stabilization of polynomial ideals over
rings~\cite{Ka97} under specialization. The published algorithm
provides a highly complex CGS, containing even incompatible
segments, but the Risa/Asir implementation makes an initial
reduction and gives a better output. We implemented an extra
routine to further reduce the output by transforming
specifications into red-specifications characterized by a pair
$(N,W)$, where $N$ is the null-condition ideal and $W$ is a set of
irreducible polynomials.

Among the tests we have done we explain four interesting ones.
\begin{example}\label{conicexss}
First we consider a very simple but illustrative example: the
discussion of the singular points of a conic already studied in
example~\ref{conicex}.

Using the Risa/Asir implementation of SACGB together with the
additional simplifications we obtain the following description of
the CGS: 
\vspace{3mm}

\begin{center}
\begin{tabular}{||c|l|l||}
 \hline \hline
 $\lpp$'s & Basis & Description  \\
 \hline
$[1]$ & $[1]$ & $\Compl^3 \setminus \V(bcd)$\\
\hline $[1]$ & $[1]$ & $\V(c) \setminus \V(d)$
\\
\hline $[x]$ & $[x+c y]$ & $\V(d,b-c^2)$
\\
 \hline $[y, x]$ & $[y, x]$ & $\V(d) \setminus \V((b-c^2)c)$
 \\ \hline $[y, x]$ & $[y, 2x+d]$ & $\V(d, c) \setminus \V(b)$
  \\ \hline $[y, x]$ & $[2c y+d, x]$ & $\V(b) \setminus \V(c d)$ \\
 \hline
\hline
\end{tabular}
\end{center}\vspace{3mm}

There are two segments with basis $[1]$, i.e. when the conic has
no singular points. The first one corresponds to the whole
$\Compl^3$ space except the three planes $\V(b)$, $\V(c)$ and
$\V(d)$. The second one corresponds to the plane $\V(c)$ except
the line $\V(c,d)$. They have empty intersection and its union
describes the unique generic segment in MCCGS, namely the whole
$\Compl^3$ space except the two planes $\V(b)$ and $\V(d)$ plus
the line $\V(b,c)$ except the origin $(0,0,0)$.

The segment with $\lpp$ set $[x]$ (i.e. the conic is a double line
of singular points) coincides with the one in MCCGS.

Finally, there are three segments with $\lpp$ set $[x,y]$, i.e.
the conic has one single singular point. The first one corresponds
to the plane $\V(d)$ minus the line $\V(c,d)$ and the parabola
$\V(d,b-c^2)$. The second one corresponds to the line $\V(c,d)$
minus the origin $(0,0,0)$. The third one corresponds to the plane
$\V(b)$ minus the lines $\V(b,c)$ and $\V(b,d)$. These three sets
have no common intersection and their union describes the plane
$\V(b)$ minus the line $\V(b,c)$ plus the plane $\V(d)$ minus the
parabola $\V(d,b-c^2)$, which is the unique segment in MCCGS. Also
the basis given by MCCGS for this segment specializes to the bases
of the three segments provided by SACGB.

Using Reduce implementation of Weispfenning's CGB, we obtained the
following CGS:\vspace{3mm}

\begin{center}
\begin{tabular}{||c|l|l||}
 \hline \hline
  Segment & Basis & Description  \\
 \hline 1 & $[bd^2]$ & $b^2cd-bc^3d\neq 0$\\
 \hline 2 & $[x^2+2cxy+dx+by^2, 2x+2cy+d,$ & \\
         & \ \ $cx+by,(2b-2c^2)y-cd]$ &
$b-c^2\neq 0, c\neq 0, bd=0$\\
 \hline 3 & $[2cdy+d^2]$ & $b\neq 0, d\neq 0, c=0$\\
 \hline 4 & $[x^2+2cxy+dx+by^2, 2x+2cy+d, cx+by]$ & $b\neq 0, c=0, d=0$ \\
 \hline 5 &  $[(2b-2c^2)y-cd]$ & $c\neq 0, d\neq 0, b-c^2=0$\\
 \hline 6 &  $[x^2+2cxy+dx+by^2, 2x+2cy+d,cx+by]$ & $c\neq 0, d\neq 0, b-c^2=0$ \\
 \hline 7 & $[4cxy+4by^2-2cdy-d^2]$ & $d\neq 0, b=0, c=0$\\
 \hline 8 & $[x^2+2cxy+dx+by^2, 2x+2cy+d]$ & $b= 0, c=0, d= 0$\\
 \hline
\hline
\end{tabular}
\end{center}\vspace{3mm}

As it can be seen, the description of the segments is not very
friendly. In order to interpret these CGS as a partition we have
manually built the following binary table in which $0$ represents
"being equal to $0$", and $1$ "being different from $0$". The last
column matches each CGS segment with one of the three MCCGS
segments identified by its $\lpp$.\vspace{3mm}

\begin{center}
\begin{tabular}{||c|c|c|c|c|c||}
 \hline \hline
  Segment & $b$ & $c$ & $d$ & $b-c^2$ & MCCGS \  lpp \\
 \hline 1 & 1 & 1 & 1 & 1 & $[1]$\\
 \hline 2 & 0 & 1 & 0 & 1 & $[x,y]$\\
          & 0 & 1 & 1 & 1 & \\
          & 1 & 1 & 0 & 1 & \\
 \hline 3 & 1 & 0 & 1 & 1 & $[1]$\\
 \hline 4 & 1 & 0 & 0 & 1 & $[x,y]$\\
 \hline 5 & 1 & 1 & 1 & 0 & $[1]$\\
 \hline 6 & 1 & 1 & 0 & 0 & $[x]$\\
 \hline 7 & 0 & 0 & 1 & 0 & $[1]$\\
 \hline 8 & 0 & 0 & 0 & 0 & $[x]$\\
 \hline
\hline
\end{tabular}
\end{center}\vspace{3mm}

The CPU times are 1.46 sec for MCCGS, 0.18 sec for SACGS and 0.05
sec for CGB.

We see that MCCGS outputs a simpler discussion, not only
theoretically but also experimentally as all the segments
corresponding to the same set of solutions are summarized in a
single segment, while SACGS and CGB do not. Nevertheless, SACGS
and CGB are both correct and faster than MCCGS, and although they
do not ensure that the $S_i$ form a partition of the parameter
space, in this example they do.
\end{example}

\begin{example}
We consider now an example proposed in~\cite{SuSa06} for which
they give the following comprehensive Gr\"obner basis wrt
$lex(t,x,y)$
\[S:=[x^3-a,y^4-b,x+y-t]\]
and ask for the CGS of $\langle S \rangle$ wrt $\lex(x,y,t)$.

MCCGS provides in 632 sec. the canonical tree shown in
Figure~\ref{SuSaEx} with only 4 segments which takes 16 lines of a
Maple worksheet.

\begin{figure}
\begin{center}
\includegraphics[width=13cm]{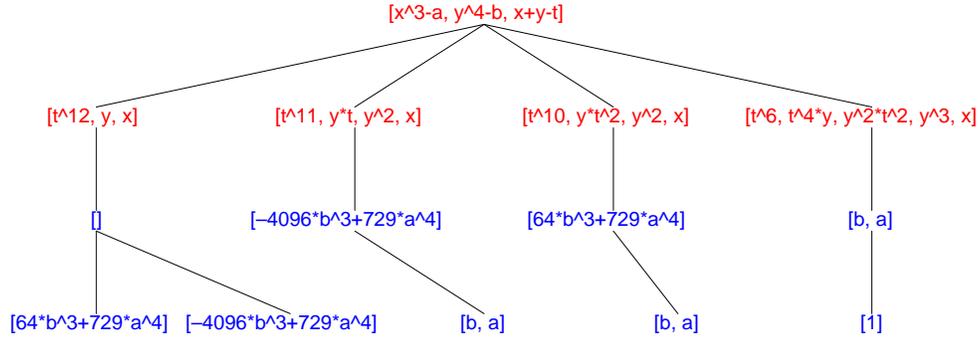} \caption{\label{SuSaEx} Canonical tree for $[x^3-a,y^4-b,x+y-t]$ wrt $\lex(x,y,t)$.}
\end{center}
\end{figure}

On the other hand, the Risa/Asir SACGB with the respective
simplifications produces in 1.62 sec. the following CGS which
takes 22 lines of a Maple worksheet: (for space restrictions we do
not print the bases)\vspace{3mm}

\begin{center}
\begin{tabular}{||c|l|l||}
 \hline \hline
 $\lpp$'s & Description  \\
 \hline $[t^{12}, y, x]$ & $\Compl^2 \setminus \V(ab(729a^4+64b^3)(729a^4-4096b^3)(16767a^4+5632b^3))$\\
 \hline $[t^{12}, y, x]$ &  $\V(16767a^3+5632b^3) \setminus \V(a b)$\\
 \hline $[t^{12}, y, x]$ &  $\V(a) \setminus \V(b)$\\
 \hline $[t^{12}, y, x]$ &  $\V(b) \setminus \V(a)$\\
 \hline $[t^{11}, ty,y^2, x]$ & $\V(729a^4-4096b^3) \setminus \V(a b)$\\
 \hline $[t^{10}, t^2y, y^2, x]$ & $\V(729a^4-64b^3) \setminus \V(a b)$\\
 \hline $[t^6, t^4y, t^2y^2, y^3, x]$ & $\V(b,a)$\\
 \hline
\hline
\end{tabular}
\end{center}\vspace{3mm}

The first segment is described by the whole $\Compl^2$ space minus
the three curves $\V(729a^4+64b^3)\text{,
}\V(729a^4-4096b^3)\text{, } \V(16767a^4+5632b^3)$ and the lines
$\V(a)$ and $\V(b)$. The second one is described by the curve
$\V(16767a^4+5632b^3)$ except the origin. The third one is the
line $\V(a)$ minus the origin and the forth segment is described
by the line $\V(b)$ minus the origin. These four segments have
empty intersection and are associated to bases with lpp set
$[t^{12}, y, x]$. Their union corresponds to the unique generic
segment in MCCGS, namely the whole $\Compl^2$ space minus the two
curves $\V(729a^4+64b^3)$ and $\V(729a^4-4096b^3)$.

The segment in SACGB with lpp set $[t^{11}, ty,y^2, x]$ is
described by the curve $\V(729a^4-4096b^3)$ except the origin,
which corresponds exactly to the segment associated to the same
lpp set in MCCGS.

The segment with lpp set $[t^{10}, t^2y, y^2, x]$ and described by
the curve $\V(729a^4-64b^3)$ minus the origin also coincides with
the one in MCCGS associated to this lpp set.

And finally, the segment having basis with lpp set $[t^6, t^4y,
t^2y^2, y^3, x]$ is described on the origin $\V(b,a)$, which
agrees with the segment associated to the same lpp set in MCCGS.

All seven segments have no common intersection and thus they form
a partition of the $\Compl^2$ space, even though SACGB does not
ensure it.
\end{example}

\begin{example}\label{exninepoints}
We also have tried to test SACGB with the systems of the nine
points circle theorem explained in section
\ref{applicationsection} above. SACGB after 3 hours of computation
went out of memory and had not yet reached an end, while MCCGS
takes only 11.45 sec. for testing the compatibility of the
hypotheses and 2.21 sec. for discussing the theorem thesis.
\end{example}

\begin{example}\label{exromin}
The last test is the system of the Romin robot\cite{GR93}:
\[ R=[a+d s_1,b-d c_1,l_2 c_2+l_3 c_3-d,l_2 s_2+l_3
s_3-c,s_1^2+c_1^2-1,s_2^2+c_2^2-1,s_3^2+c_3^2-1] \] wrt
$\lex(c_3,s_3,c_2,s_2,c_1,s_1)$ and $\lex(l_2,l_3,a,b,c,d)$. MCCGS
takes 43.23 sec in discussing the system and provides 9 segments.
SACGB also went out of memory.
\end{example}

\begin{section}*{Conclusions}
The interest of MCCGS relies, essentially, in the simplicity of
the output for applications, and in the canonical character of it,
conceding an easier interpretation of the results. We have also
observed that the obtention of the MCCGS from the BUILDTREE CGS
only increases the computation time in about 20-30\%.

The existence of the MCCGS depends on the Conjecture formulated
in~\cite{Mo07}. The use of the algorithm will provide evidence of
it or a counterexample. In almost all the high number of tests
that we have done the algorithm has always obtained a unique
segment for each different $\lpp$ set, confirming the conjecture.
The only ideal for which the algorithm obtains two different
segments with the same $\lpp$ is $=\langle u(ux+1), (ux+1)x
\rangle$ proposed by~\cite{Wi06}, and there both segments are
clearly intrinsically different and cannot be merged nor
summarized into a single one. Thus this example also provides
evidence of the Conjecture. To give a counterexample proving the
falsehood of the Conjecture, we must find an ideal for which the
algorithm MCCGS obtains two or more segments with the same $\lpp$
which could be merged or summarized in a different way.

Although we have only made some very simple tests with CGB, we
have observed that it seems faster than SACGB and MCCGS in those
specific problems. It stands out for computing a CGS with faithful
bases which are not always useful for applications.
Experimentally, it seems to obtain a partition of the parameter
space, even if there is no theoretic evidence. Nevertheless, the
number of segments is much higher than MCCGS and are difficult to
understand.

SACGB stands out for being in general very reliable to compute a
CGS. Its efficiency depends on the type of system to be dealt
with. It seems to behave faster than MCCGS in problems for which a
low number of cases is expected. Furthermore, we must remind that
the output of SACGB is very complex and also needs extra
simplifications to be interpreted.

One can also adapt the MCCGS algorithm to the CGS obtained by
other algorithms instead of BUILDTREE. To do this one needs to
transform the output of the involved algorithm into a disjoint
reduced CGS, and then apply step ii) and iii), i.e. SELECTCASES
and MCCGS.

MCCGS takes, generally, more CPU time for simple problems.
Nevertheless the simplifications inside MCCGS often allow to
discuss systems of higher complexity, as seen in examples
\ref{exninepoints} and \ref{exromin} above.

Finally, we have seen that MCCGS algorithm stands out for having
the best features to be used for automatic theorem proving and
discovering as well as for other applications.
\end{section}

\section*{Ackowledgements} We are very grateful to Josep M.
Brunat for his suggestions on the recursive handling of trees. We
would also like to thank Professor Pelegr\'{\i} Viader for his
many helpful comments and his insightful perusal of our first
draft. Finally, we also want to thank our referees for their very
helpful comments.


\begin{thebibliography}{99}


   \bibitem[Dolzman-Seidl-Sturm 2006]{DoSeSt06} A. Dolzmann, A. Seidl, T. Sturm. (2006)
   \newblock REDLOG software in REDUCE http://staff.fim.uni-passau.de/$\sim$\,sturm/


   \bibitem[Gonz\'alez-Recio 1993]{GR93} Gonz\'{a}lez-L\'{o}pez, M.J.,  Recio,
   T., 1993.
   \newblock The ROMIN inerse geometric model and the dynamic evaluation
   method.
   \newblock In: Computer Algebra in Industry, A.M. Cohen ed., John Wiley \& Sons: 117--141 (1993).

   \bibitem[Gonz\'alez-Traverso-Zanoni 2005]{GoTrZa05} Gonz\'alez-Vega, L., Traverso, C., Zanoni,
   A., 2005.
   \newblock Hilbert Stratification and Parametric Groebner Bases.
   \newblock CASC-2005, p. 220-235.

   \bibitem[Kalkbrenner 1997]{Ka97} Kalkbrenner, M., 1997.
   \newblock On the stability of Gr\"obner bases under
   specializations.
   \newblock {\em Jour. Symb. Comp.}, {\bf 24}(1): (1997), 51--58.

   \bibitem[Manubens-Montes 2006]{MaMo06} Manubens, M., Montes, A., 2006.
   \newblock Improving DISPGB Algorithm Using the Discriminant Ideal.
   \newblock {\em Jour. Symb. Comp.} {\bf 41} (2006), 1245--1263.

   \bibitem[Montes 2002]{Mo02} Montes, A., 2002.
   \newblock New Algorithm for Discussing Gr\"obner Bases with Parameters.
   \newblock {\em Jour. Symb. Comp.} {\bf 33}:1-2 (2002), 183--208.

   \bibitem[Montes 2007]{Mo07} Montes, A., 2007.
   \newblock On the canonical discussion of polynomial systems with parameters.
   \newblock Preprint arXiv: AC/0601674.

   \bibitem[Montes-Recio 2007]{MoRe07} Montes, A., Recio, T., 2006.
   \newblock Automatic discovery of geometry theorems using minimal canonical comprehensive Groebner systems.
   \newblock Preprint arXiv: AC/0703483.

   \bibitem[Nabeshima 2005]{Na05} Nabeshima, K., 2005.
   \newblock A computation method for ACGB-V.
   \newblock Proceedings of A3L 2005 (Conference in Honour of the 60th
   Bithday of V. Weispfenning), eds. A. Dolzmann, A Seidl, T.
   Sturm. p 171-180. BOD Norderstedt.

   \bibitem[Nabeshima 2006]{Na06} Nabeshima, K., 2006.
   \newblock Comprehensive Groebner Bases for Modules. Short communication, ACA-2006, Varna.

   \bibitem[Sato 2005]{Sa05} Sato, Y., 2005.
   \newblock Stability of Gr\"obner basis and ACGB.
   \newblock Proceedings of A3L 2005 (Conference in Honour of the 60th
   Birthday of V. Weispfenning), eds. A. Dolzmann, A Seidl, T.
   Sturm. p 223-228. BOD Norderstedt.

   \bibitem[Sato-Suzuki 2003]{SaSu03} Sato, Y., Suzuki, A., 2003.
   \newblock An alternative approach to Comprehensive Gr\"obner bases.
   \newblock {\em Jour. Symb. Comp.} {\bf 36}:3-4 (2003), 649-667.

   \bibitem[Suzuki-Sato 2006]{SuSa06} Suzuki, A., Sato, Y., 2006.
   \newblock A Simple Algorithm to compute Comprehensive Gr\"obner
   bases. Proceedings of ISSAC 2006, ACM. p 326-331.
   \newblock Implementation in Risa/Asir and Maple
   (http://kurt.scitec.kobe-u.ac.jp/$\sim$sakira/).

   \bibitem[Weispfenning 1992]{We92} Weispfenning, V., 1992.
   \newblock Comprehensive Gr\"obner bases.
   \newblock {\em Jour. Symb. Comp.} {\bf 14}:1-1 (1992), 1-29.

   \bibitem[Weispfenning 2003]{We03} Weispfenning, V., 2003.
   \newblock Canonical Comprehensive Gr\"obner bases.
   \newblock {\em Jour. Symb. Comp.} {\bf 36}:3-4 (2003), 669-683.

   \bibitem[Wibmer 2006]{Wi06} M.~Wibmer.
   \newblock Gr\"obner bases for families of affine schemes.
   \newblock arXiv. math/0608019 (2006).

\end{thebibliography}
\end{document}